\documentclass{elsart3-1}

\usepackage{amssymb}

\usepackage[english,francais]{babel}

\newtheorem{theorem}{Theorem}[section]

\newtheorem{e-proposition}[theorem]{Proposition}

\newtheorem{e-definition}[theorem]{Definition\rm}

\newtheorem{theoreme}{Th{\'e}or{\`e}me}[section]

\newtheorem{proposition}[theoreme]{Proposition}

\renewcommand{\theequation}{\arabic{equation}}
\def\og{\leavevmode\raise.3ex\hbox{$\scriptscriptstyle\langle\!\langle$~}}
\def\fg{\leavevmode\raise.3ex\hbox{~$\!\scriptscriptstyle\,\rangle\!\rangle$}}

\begin{document}

\begin{frontmatter}

% Titre, auteurs et adresses

% utiliser la commande \thanksref dans \title, \author ou \address
%     pour les notes en bas de page ;

% utiliser la commande \ead pour l'adresse e-mail de chaque auteur
%    (apr{\"E}s la commande \auteur) ;

% \author{Name\thanksref{label2}}
\selectlanguage{francais}
\title{Th{\'e}or{\`e}me d'{\'e}quidistribution de Brolin en dynamique $p$-adique}

%\vspace{-2.6cm}

\selectlanguage{english}
\title{Brolin's equidistribution theorem in $p$-adic dynamics}

% utiliser les {\`E}tiquettes pour indiquer l'adresse de chaque auteur,
%     s'il y a plusieurs adresses


\author[CF]{Favre Charles}
\ead{favre@math.jussieu.fr}
\author[JRL]{Juan Rivera-Letelier\thanksref{soutien}}
\thanks[soutien]{partiellement soutenu par le projet FONDECYT N 1040683}
\ead{juanrive@ucn.cl}

\address[CF]{
CNRS et Institut de Math{\'e}matiques de Jussieu,
Case 7012,
2  place Jussieu,
F-75251 Paris Cedex 05 - France}
\address[JRL]{Departamento de Matem{\'a}tica,
Universidad Cat{\'o}lica del Norte,
Casilla 1280,
Antofagasta - Chile} 

% etc, etc

\begin{abstract}
% Texte de l'abstract en anglais
  We prove an analog of the famous equidistribution theorem of Brolin
  for rational mappings in one variable defined over the $p$-adic
  field $\mathbb{C}_p$. We construct a mixing invariant probability
  measure which describes the asymptotic distribution of iterated preimages of
  a given point. This measure is supported on the Berkovich space
  $\mathsf{P}^1(\mathbb{C}_p)$ associated to
  $\mathbb{P}^1(\mathbb{C}_p)$. We show that its support is precisely
  the Julia set of $R$ as defined in~\cite{R3}.  Our results are based
  on the construction of a Laplace operator on real trees with
  arbitrary number of branching as done in~\cite{FJ}.

\vskip 0.5\baselineskip

\selectlanguage{francais}
% Texte du r{\`E}sum{\`E} en fran{\'A}ais
\noindent{\bf R{\'e}sum{\'e}}
\vskip 0.5\baselineskip
\noindent
Nous d{\'e}montrons un analogue du th{\'e}or{\`e}me classique
d'{\'e}quidistribution de Brolin pour les applications rationnelles {\`a}
une variable d{\'e}finies sur le corps $p$-adique $\mathbb{C}_p$.  On
construit une mesure invariante et m{\'e}langeante qui d{\'e}crit la
distribution (asymptotique) des pr{\'e}images it{\'e}r{\'e}es d'un point
donn{\'e}.  Cette mesure est {\`a} support dans l'espace analytique de
$\mathbb{P}^1(\mathbb{C}_p)$, au sens de Berkovich, que l'on note
$\mathsf{P}^1(\mathbb{C}_p)$.  On d{\'e}montre que le support de cette
mesure est {\'e}gale {\`a} l'ensemble de Julia dans
$\mathsf{P}^1(\mathbb{C}_p)$, introduit dans~\cite{R3}.  Nos
r{\'e}sultats sont bas{\'e}s sur la notion d'op{\'e}rateur de Laplace sur
les arbres r{\'e}els avec  nombre arbitraire de branchements construit
dans~\cite{FJ}.

\end{abstract}
\end{frontmatter}

% Maintenant la version abr{\`E}g{\`E}e en anglais, si pr{\`E}sente
\selectlanguage{english}
\section*{Abridged English version}
% Texte de la version abr{\`E}g{\`E}e en anglais
Iteration theory of rational maps with coefficients in $\mathbb{C}_p$
has recently received increasing attention, see for instance
~\cite{Be1},~\cite{Be2},~\cite{BH},~\cite{Hs},~\cite{HY},
~\cite{MS},~\cite{PST},~\cite{R1},~\cite{R2},~\cite{Y}. The local
dynamics is now well-understood, and the action of a rational map $R$
on its Fatou set has been studied in great details. Ergodic properties
(topological entropy, existence of mixing measures) have raised
however much less interest, except in connection with number theory,
see~\cite{BH},~\cite{PST}.  In the complex case some (mixing)
invariant measures are constructed from potential theoretic
considerations, but such a theory is not available on $\mathbb{C}_p$,
see however~\cite{Ru}.  An essential problem lies in the fact that
$\mathbb{C}_p$ endowed with its $p$-adic norm is not locally compact,
and is totally discontinuous. In order to remedy to these problems, we
follow the approach of Berkovich, and add points to $\mathbb{C}_p$ in
order to obtain an arcwise connected compact topological space suited
for analysis.  More precisely, one considers the set of valuations
$\nu : \mathbb{C}_p [t] \to \overline{\mathbb{R}}= \mathbb{R} \cup \{
+ \infty\}$, such that $\nu|_{\mathbb{C}_p} = \nu_p = -\log_p | \cdot
|_p$. The plane $\mathbb{C}_p$ can be naturally embedded in this space
through the map $v_x(f) = -\log_p | P(x)|_p$ for $P \in
\mathbb{C}_p[t]$.  We shall also identify the point $\infty$ in the
standard projective space $\mathbb{P}^1(\mathbb{C}_p)$ with the
function identically $+\infty$ on $\mathbb{C}_p [t]\setminus
\mathbb{C}_p^*$. With the topology of pointwise convergence (we shall
call it the weak topology in the sequel), this space of valuations
becomes a compact space we denote by $\mathsf{P}^1(\mathbb{C}_p)$. Its
topological structure has been described in great details
in~\cite{Ber},~\cite{R2}: $\mathsf{P}^1(\mathbb{C}_p)$ is a \emph{non
  metric real tree} in the sense of~\cite[chapitre 3]{FJ} (or a simply
connected quasi-polyhedron according to~\cite{Ber}). Points of
$\mathbb{P}^1(\mathbb{C}_p)$ are end points in the tree, and we denote
the ``interior'' of the tree by $\mathbb{H}_p =
\mathsf{P}^1(\mathbb{C}_p)\setminus \mathbb{P}^1(\mathbb{C}_p)$. This
set is again a real tree, which posesses a natural metric $d$. This
metric is uniquely characterized by the fact that it is invariant by
the action of $\mathrm{PGL}(2,\mathbb{C}_p)$, and for all valuations
$\nu,\nu'$ on the segment $]0, + \infty[ \subset \mathbb{H}_p$, we
have $d(\nu, \nu') = |\nu(t) - \nu'(t)|$.  
%From this characterization,
%one sees that $\mathbb{H}_p$ can be also seen as the Bruhat-Tits tree
%attached to the group $\mathrm{PGL}(2,\mathbb{C}_p)$
%(see~\cite{Ber},~\cite{R2}).

Any rational map $R$ with coefficients in $\mathbb{C}_p$ admits a
natural action on $\mathsf{P}^1(\mathbb{C}_p)$ by duality: $R_* \nu
(P) = \nu ( P \circ R)$ for all $P \in \mathbb{C}_p[t]$. This action
is continuous on $\mathsf{P}^1(\mathbb{C}_p)$, and we can extend the
notions of Fatou/Julia sets to $R_*$, see~\cite{R3}.  A point $\nu$
lies in the Julia set of $R_*$ if and only if for all weak open set
$U\subset \mathsf{P}^1(\mathbb{C}_p)$ containing $\nu$, the set
$\cup_{n\ge0} R^n(U)$ contains $\mathbb{H}_p$. It is also the closure
of the set of repelling periodic points (in a suitable sense for
points in $\mathbb{H}_p$, see~\cite{R3}).

If $R$ is a rational map and $z\in \mathbb{P}^1(\mathbb{C}_p)$, we
denote by $\delta_z$ the Dirac mass at $\{z\}$, and we define $R^*
\delta_z $ as the atomic measure supported on $R^{-1}\{ z \}$ whose
mass at a point $w \in R^{-1}\{ z \}$ is equal to the local
topological degree of $R$ at $w$. The total mass of $R^* \delta_z$ is
always equal to the degree $D$ of  $R$.

Our main result is the following
\begin{theorem}
  For any rational map $R$ of degree $D\ge 2$ and defined over
  $\mathbb{C}_p$, there exists a probability measure  $\rho_R$ on
  $\mathsf{P}^1(\mathbb{C}_p)$, which is invariant by $R_*$, mixing, whose
  support equals the Julia set of $R_*$ in
  $\mathsf{P}^1(\mathbb{C}_p)$, and such that  $ \lim_{n\to\infty}
  D^{-n} R^{n*}\delta_z= \rho_R~, $ for each point $z\in
  \mathbb{P}^1(\mathbb{C}_p)$ which is not totally invariant by 
  $R$ or $R^2$.
\end{theorem}

Although the measure $\rho_R$ appears for the first time in the
present work, a weaker form of our distribution result has been
formulated and proved in~\cite{PST}.

Note that the set $\mathcal{E}$ of points of
$\mathbb{P}^1(\mathbb{C}_p)$ which are totally invariant by $R$ or
$R^2$ consists of at most two points.  When $\mathcal{E}$ has two
points, then $R$ or is conjugated by an element in $\mathrm{PGL}(2,
\mathbb{C}_p)$ to $z \mapsto z^{\pm D}$; when $\mathcal{E}$ is a singleton,
$R$ is conjugated to a polynomial map. One can also prove a stronger
version of the theorem, namely $ \lim_{n\to\infty} D^{-n} R^{n*}\rho=
\rho_R$ for any probability measure $\rho$ which does not charge
$\mathcal{E}$.

Following~\cite{MS}, we say that $R$ has good reduction when its
reduction in the residue field of $\mathbb{C}_p$ (i.e. the algebraic
closure of $\mathbb{Z}/ p \mathbb{Z}$) is well-defined and has the
same degree as $R$.  When it is the case, the valuation $\nu_0$
determined by the conditions $\nu_0( t - z) = \min \{ 0, \nu_p(z) \}$
for all $z \in \mathbb{C}_p$ is totally invariant by $R_*$, so $
\rho_R = \delta_{\nu_0}$. Note that the Julia set of $R_*$ is then
reduced to $\{ \nu_0 \}$, see~\cite{R3}.  This situation can be
generalized to the case when $R$ is conjugated by some element in
$\mathrm{PGL}(2, \mathbb{C}_p)$ to a rational map having good
reduction.
\begin{proposition}
  The measure $\rho_R$ does not charge points of
  $\mathbb{P}^1(\mathbb{C}_p)$.  Moreover, $\rho_R\{ \nu \} >0$ for
  some $\nu \in \mathbb{H}_p$, if and only if one can find $\phi \in
  \mathrm{PGL}(2, \mathbb{C}_p)$ such that $\phi \circ R \circ
  \phi^{-1}$ has good reduction. In this case, $\phi (\nu ) = \nu_0$
  and $\rho_R = \delta_{\nu}$.
\end{proposition}

Details of proofs, more precise ergodic properties of $\rho_R$,
applications to the computation of normalized heights in the spirit
of~\cite{PST}, and to the distribution of points of small height will
appear in a later work.

\selectlanguage{francais}
% texte principale
\section{Introduction}
\label{Sintro}

L'it{\'e}ration des fractions rationnelles {\`a} coefficients dans le corps
$p$-adique $\mathbb{C}_p$ conna{\^\i}t depuis quelques temps un int{\'e}r{\^e}t
croissant (voir par
exemple~\cite{Be1},~\cite{Be2}, \cite{BH},~\cite{Hs},~\cite{HY},
~\cite{MS},~\cite{PST},~\cite{R1},~\cite{R2},~\cite{Y}). La dynamique
locale au voisinage d'un point fixe est maintenant compl{\`e}tement
comprise, et le comportement d'une fraction rationnelle $R$ sur son
domaine de Fatou a {\'e}t{\'e} {\'e}tudi{\'e} en d{\'e}tail. Les propri{\'e}t{\'e}s ergodiques de
$R$ (entropie topologique, existence de mesures m{\'e}langeantes) ont
cependant attir{\'e}es peu d'attention.  Dans le cas complexe, on utilise
g{\'e}n{\'e}ralement des m{\'e}thodes d'analyse harmonique pour construire des
mesures invariantes, mais une telle th{\'e}orie est tr{\`e}s peu d{\'e}velopp{\'e}e
sur $\mathbb{C}_p$ (voir cependant~\cite{Ru}). Un probl{\`e}me essentiel
provient du fait que, muni de sa norme $p$-adique $| \cdot |_p$, le
corps $\mathbb{C}_p$ n'est pas un espace localement compact, et est de
surcro{\^\i}t totalement discontinu. Pour rem{\'e}dier {\`a} cel{\`a}, suivant
Berkovich~\cite{Ber}, on rajoute {\`a} $\mathbb{C}_p$ des points afin de
le \og connexifier \fg~, ce qui permet de d{\'e}velopper des outils
d'analyse sur $\mathbb{C}_p$ analogues {\`a} ceux existant sur
$\mathbb{C}$.  Plus pr{\'e}cis{\'e}ment, on consid{\`e}re l'espace des valuations
$\nu : \mathbb{C}_p [t] \to \overline{\mathbb{R}}= \mathbb{R} \cup \{
+ \infty\}$, telles que $\nu|_{\mathbb{C}_p} = \nu_p = -\log_p | \cdot
|_p$. Le plan $\mathbb{C}_p$ se plonge alors naturellement dans cet
espace par le morphisme d'{\'e}valuation $v_x(P) = -\log_p | P(x)|_p$ pour
$P \in \mathbb{C}_p[t]$.  On identifiera aussi le point $\infty$ de
l'espace projectif standard $\mathbb{P}^1(\mathbb{C}_p)$ {\`a} la fonction
identiquement $+\infty$ sur $\mathbb{C}_p [t]\setminus
\mathbb{C}_p^*$. Muni de la topologie de la convergence simple (on
dira aussi topologie faible dans la suite), cet ensemble de valuations
devient un espace compact que l'on notera
$\mathsf{P}^1(\mathbb{C}_p)$. Sa structure topologique a {\'e}t{\'e} d{\'e}crite
en grand d{\'e}tail dans~\cite{Ber},~\cite{R2}:
$\mathsf{P}^1(\mathbb{C}_p)$ est un \emph{arbre r{\'e}el non m{\'e}trique} au
sens de~\cite[chapitre 3]{FJ} (ces espaces sont d{\'e}nomm{\'e}s aussi
quasi-polyh{\`e}dre simplement connexe dans~\cite{Ber}). Les points de
$\mathbb{P}^1(\mathbb{C}_p)$ en constituent le bord (ce sont les bouts
de l'arbre), et on notera son \og int{\'e}rieur\fg~ $\mathbb{H}_p =
\mathsf{P}^1(\mathbb{C}_p)\setminus \mathbb{P}^1(\mathbb{C}_p)$. Cet
ensemble est {\`a} nouveau un arbre r{\'e}el, et il poss{\`e}de une m{\'e}trique
naturelle que l'on notera $d$. Elle est caract{\'e}ris{\'e}e par le fait
d'{\^e}tre invariante par l'action de $\mathrm{PGL}(2,\mathbb{C}_p)$, et
pour toutes valuations $\nu,\nu'$ sur le segment $]0, + \infty[
\subset \mathbb{H}_p$, on a $d(\nu, \nu') = |\nu(t) - \nu'(t)|$.  
%De cette caract{\'e}risation, on voit que $\mathbb{H}_p$ s'interpr{\`e}te
%aussi comme l'arbre de Bruhat-Tits associ{\'e} au groupe
%$\mathrm{PGL}(2,\mathbb{C}_p)$ (voir~\cite{Ber},~\cite{R2}).

Toute application rationnelle $R$ {\`a} coefficients dans $\mathbb{C}_p$
admet une action naturelle sur $\mathsf{P}^1(\mathbb{C}_p)$ par
dualit{\'e}: $R_* \nu (P) = \nu ( P \circ R)$ pour tout $P \in
\mathbb{C}_p[t]$. Cette action est continue sur
$\mathsf{P}^1(\mathbb{C}_p)$, et l'on peut {\'e}tendre les notions
d'ensembles de Fatou/Julia pour $R_*$, voir~\cite{R3}.  Un point $\nu$
est dans l'ensemble de Julia de $R_*$ si et seulement si pour tout
ouvert faible $U\subset \mathsf{P}^1(\mathbb{C}_p)$ contenant $\nu$,
l'ensemble $\cup_{n\ge0} R^n(U)$ contient $\mathbb{H}_p$. C'est aussi
la clot{\^u}re de l'ensemble des points p{\'e}riodiques r{\'e}pulsifs (en un sens
convenable pour les points de $\mathbb{H}_p$, voir~\cite{R3}).

Si $R$ est une application rationnelle, et $z\in
\mathbb{P}^1(\mathbb{C}_p)$, on note $\delta_z$ la masse de Dirac
support{\'e}e sur $\{z\}$, et on d{\'e}finit $R^* \delta_z $ comme la
mesure atomique support{\'e}e sur $R^{-1}\{ z \}$ dont la masse en un
point $w \in R^{-1}\{ z \}$
co{\"\i}ncide avec le degr{\'e} topologique local de $R$ en $w$. La
masse de $R^* \delta_z$ est toujours {\'e}gale au degr{\'e} $D$ de $R$.

Notre  r{\'e}sultat principal est le suivant.
\begin{theoreme}
  Pour toute fraction rationnelle $R$ de degr{\'e} $D\ge 2$ et d{\'e}finie sur
  $\mathbb{C}_p$, il existe une mesure de probabilit{\'e} $\rho_R$ sur
  $\mathsf{P}^1(\mathbb{C}_p)$, invariante par $R_*$ et m{\'e}langeante,
  dont le support co{\"\i}ncide avec l'ensemble de Julia de $R_*$ dans
  $\mathsf{P}^1(\mathbb{C}_p)$, et telle que $ \lim_{n\to\infty}
  D^{-n} R^{n*}\delta_z= \rho_R~, $ pour tout point $z\in
  \mathbb{P}^1(\mathbb{C}_p)$ qui n'est pas totalement invariant par
  $R$ ou $R^2$.
\end{theoreme}
Bien que la mesure $\rho_R$ elle-m{\^e}me apparaisse pour la
premi{\`e}re fois dans ce travail, il convient de noter qu'une forme
plus faible du r{\'e}sultat d'{\'e}quidistribution ci-dessus a
{\'e}t{\'e} formul{\'e} dans~\cite[section 5.4]{PST}.

Notons que l'ensemble $\mathcal{E}$ des points de
$\mathbb{P}^1(\mathbb{C}_p)$ totalement invariants par $R$ ou $R^2$
consiste d'au plus deux points.  Lorsque $\mathcal{E}$ en poss{\`e}de
deux, $R$ est conjugu{\'e} par un {\'e}l{\'e}ment de $\mathrm{PGL}(2,
\mathbb{C}_p)$ {\`a} $z \mapsto z^{\pm D}$; lorsque $\mathcal{E}$ en
poss{\`e}de un, $R$ est conjugu{\'e} {\`a} un polyn{\^o}me. On peut de
plus d{\'e}montrer que $ \lim_{n\to\infty} D^{-n} R^{n*}\rho= \rho_R$
pour toute mesure de probabilit{\'e} $\rho$ ne chargeant pas
$\mathcal{E}$.

Suivant~\cite{MS}, on dit que $R$ a bonne r{\'e}duction si et
seulement si sa r{\'e}duction dans le corps r{\'e}siduel de
$\mathbb{C}_p$ (i.e. la cl{\^o}ture alg{\'e}brique de $\mathbb{Z}/ p
\mathbb{Z}$) est bien d{\'e}finie et de m{\^e}me degr{\'e} que $R$.
Lorsque c'est le cas, la valuation $\nu_0$ d{\'e}termin{\'e}e par
$\nu_0( t - z) = \min \{ 0, \nu_p(z) \}$ pour tout $z \in
\mathbb{C}_p$ est totalement invariante par $R$, et donc $ \rho_R =
\delta_{\nu_0}$.  Notons que l'ensemble de Julia de $R$ est alors
r{\'e}duit au point $\{ \nu_0 \}$, voir~\cite{R3}.  Cette situation se
g{\'e}n{\'e}ralise imm{\'e}diatement au cas o{\`u} $R$ est
conjugu{\'e}e par un {\'e}l{\'e}ment de $\mathrm{PGL}(2,
\mathbb{C}_p)$ {\`a} une fraction rationnelle ayant bonne
r{\'e}duction.
\begin{proposition}\label{prop}
  La mesure $\rho_R$ ne charge pas les points de
  $\mathbb{P}^1(\mathbb{C}_p)$.  De plus, $\rho_R\{ \nu \} >0$ pour un
  point $\nu \in \mathbb{H}_p$, si et seulement si on peut trouver
  $\phi \in \mathrm{PGL}(2, \mathbb{C}_p)$ tel que $\phi \circ R \circ
  \phi^{-1}$ a bonne r{\'e}duction. Dans ce cas $\phi (\nu ) = \nu_0$ et
  $\rho_R = \delta_\nu$.
\end{proposition}
Les d{\'e}tails des preuves, ainsi que les propri{\'e}t{\'e}s
ergodiques de $\rho_R$, les applications au calcul de hauteurs
normalis{\'e}es dans l'esprit de~\cite{PST}, et \`a l'\'equidistribution
des points de petites hauteurs para{\^\i}tront dans un travail
ult{\'e}rieur.

\section{Op{\'e}rateur de Laplace}
Nous allons donner les grandes lignes de la construction d'une classe
de fonctions $\mathcal{P}$ d{\'e}finies sur $\mathbb{H}_p$ {\`a}
valeurs r{\'e}elles, et d'un op{\'e}rateur not{\'e} $\Delta$ qui {\`a}
chaque fonction $g$ de $\mathcal{P}$ associe une mesure bor{\'e}lienne
sign{\'e}e d{\'e}finie sur $\mathsf{P}^1(\mathbb{C}_p)$,
diff{\'e}rence de deux mesures positives de masse {\'e}gales.
L'op{\'e}rateur $\Delta$ est appel{\'e} \emph{op{\'e}rateur de
Laplace}, et est un analogue de l'op{\'e}rateur de Laplace standard
sur un graphe fini ou sur $\mathbb{R}^n$, adapt{\'e} ici {\`a} la
g{\'e}om{\'e}trie d'arbre de $\mathbb{H}_p$. Pour les arbres finis, un
tel op{\'e}rateur se construit facilement; mais dans notre contexte
les points de branchements sont denses sur tout segment, et le nombre
de branches est infini en chacun de ces points.  La m{\'e}thode de
construction suit tr{\`e}s pr{\'e}cis{\'e}ment~\cite[Chapitre 7]{FJ}.

Fixons $\nu_*\in\mathbb{H}_p$. Un tel point marqu{\'e} induit une relation
d'ordre partiel $\le$ sur $\mathbb{H}_p$: $\nu \le \nu'$ d{\`e}s que
$[\nu_*,\nu]\subset [\nu_*,\nu']$. Il en est l'unique {\'e}l{\'e}ment minimal.
Dans ce contexte, on peut d{\'e}finir une notion de fonction {\`a} variations
born{\'e}es sur $(\mathbb{H}_p,\le)$, voir~\cite[Definition 7.20]{FJ}. On
{\'e}crit alors $g\in \mathcal{P}$ si l'on a $g(\nu) = \mathrm{Cst} +
\int_0^{d(\nu_*,\nu)} f(\nu_t)~dt$ avec $f$ {\`a} variations born{\'e}es, et
o{\`u} $\nu_t$ d{\'e}note l'unique valuation de $[\nu_*, \nu]$ {\`a} distance $t$
de $\nu_*$.  On v{\'e}rifie que la classe $\mathcal{P}$ ne d{\'e}pend
pas du choix de $\nu_*$.

Pour toute fonction {\`a} variations born{\'e}es $f$ suffisament r{\'e}guli{\`e}re, il
existe une unique mesure bor{\'e}lienne $\rho_f$ sur
$\mathsf{P}^1(\mathbb{C}_p)$ telle que $\rho_f\{ \mu \ge \nu \} =
f(\nu)$ par~\cite[Theorem 7.36]{FJ}. Pour $g \in \mathcal{P}$, et en
utilisant les notations pr{\'e}c{\'e}dentes on pose alors $\Delta g = 
f(\nu) \cdot \delta_\nu - \rho_f$. Cette mesure ne d{\'e}pend pas non plus du choix
de $\nu_*$.

On montre facilement que $\Delta g$ est une mesure r{\'e}elle
diff{\'e}rence de deux mesures positives de m{\^e}me masse, et que
r{\'e}ciproquement toute diff{\'e}rence de mesures positives de m{\^e}me
masse peut s'{\'e}crire sous cette forme. Si par exemple
$\nu\in\mathbb{H}_p$ et $\nu'\in\mathsf{P}^1(\mathbb{C}_p)$, et $g$
est la fonction localement constante hors de $[\nu,\nu']$ et valant
$d(\cdot,\nu)$ sur ce segment, alors $\Delta g = \delta_{\nu}
-\delta_{\nu'}$. 

Pour tout $g\in\mathcal{P}$, la masse de $\Delta g$ en un point
$\nu\in\mathbb{H}_p$ est {\'e}gale {\`a} la somme des d{\'e}riv{\'e}es de $g$ sur
toutes les branches partant de $\nu$. En un point de
$\mathbb{P}^1(\mathbb{C}_p)$, cette masse est donn{\'e}e par la limite de
la d{\'e}riv{\'e}e de $g$ quand on converge vers $\nu$ le long d'un segment.

On d{\'e}duit de~\cite[Theorem 7.61]{FJ}, la propri{\'e}t{\'e} de continuit{\'e}
suivante. Fixons $\nu\in\mathbb{H}_p$. Soit $\rho_n$ une suite de
mesures de \emph{probabilit{\'e}}, et fixons $g_n\in \mathcal{P}$ telles
que $\rho_n = \delta_\nu +\Delta g_n$. Si $g_n$ converge
ponctuellement vers $g$ dans $\mathbb{H}_p$, alors $g\in\mathcal{P}$,
et $\rho_n$ converge vaguement vers $ \delta_\nu +\Delta g$.

\section{Construction de la mesure}
Fixons une application rationnelle $R$ de degr{\'e} $D\ge2$ {\`a} coefficients
dans $\mathbb{C}_p$.  Son action sur $\mathsf{P}^1(\mathbb{C}_p)$,
not{\'e}e $R_*$, est faiblement continue, et chaque point de
$\mathsf{P}^1(\mathbb{C}_p)$ poss{\`e}de au plus $D$ pr{\'e}images.  En un
point $x\in \mathbb{P}^1(\mathbb{C}_p)$, on d{\'e}finit le degr{\'e} local,
not{\'e} $\mathrm{deg}\, _R(x)$, comme l'ordre d'annulation de la d{\'e}riv{\'e}e
de $R$ plus une unit{\'e}. Cette fonction s'{\'e}tend naturellement {\`a}
$\mathsf{P}^1(\mathbb{C}_p)$ de telle sorte que $\sum_{R_* \mu = \nu}
\mathrm{deg}\, _R (\mu) =D$, voir~\cite{R3}.

Soit $f: \mathsf{P}^1(\mathbb{C}_p)\to\mathbb{R}$ une fonction
continue.  On d{\'e}finit $R_* f( \nu) = \sum_{R_* \mu = \nu}
\mathrm{deg}\, _R (\mu) f ( \mu)$. Cette fonction est encore continue,
voir~\cite{R3}.  Si $\rho$ est une mesure de masse finie sur
$\mathsf{P}^1(\mathbb{C}_p)$, on peut alors d{\'e}finir $R^* \rho $ par
dualit{\'e} en posant $\int f\, d(R^* \rho) = \int R_*f\, d\rho$. On
v{\'e}rifie que pour tout point de $x \in \mathbb{P}^1(\mathbb{C}_p)$, on
a $R^* \delta_x = \sum_{R (y) = x} \mathrm{deg}\, _R (y) \delta_y$.
Par continuit{\'e}, on en d{\'e}duit que la masse de $R^* \rho$ est $D$ fois
la masse de $\rho$ pour toute mesure positive.  Si $x \ne x'\in
\mathbb{P}^1(\mathbb{C}_p)$, et $g$ est une fonction telle que $\Delta
g = \delta_x - \delta_{x'}$, il est facile de voir que $R^*\Delta g =
\Delta (g \circ R_*)$. On v{\'e}rifie que cette formule reste valide pour
n'importe quelle fonction $g$ dont le laplacien est une mesure
atomique, puis par continuit{\'e} pour toute fonction $g$ dans le domaine
de $\Delta$.

Prenons maintenant $\nu$ un point arbitraire de $\mathbb{H}_p$.  Nous
pouvons {\'e}crire $D^{-1}R^* \delta_\nu = \delta_\nu + \Delta g$ avec $g:
\mathbb{H}_p \to \mathbb{R}$.  Le support de $R^*\delta_\nu$ est constitu{\'e}
d'au plus $D$ points, donc l'enveloppe convexe de
l'ensemble $\mathrm{supp}\, R^*\delta_\nu \cup \{ \nu\}$ est un arbre fini
$\mathcal{T}$ dont les bouts sont situ{\'e}s dans $\mathbb{H}_p$.  Le
potentiel $g$ est localement constant en dehors de $\mathcal{T}$ et
est born{\'e} sur $\mathcal{T}$, il est donc born{\'e} sur $\mathbb{H}_p$ tout
entier.  On obtient donc
$$
\frac{1}{D^n} R^{n*}\delta_\nu = \delta_\nu + \Delta g_n 
\mbox{ avec } g_n = \sum_0^{n-1}  D^{-k} g\circ R^k_*~.
$$
La suite $g_n$ converge uniform{\'e}ment sur $\mathbb{H}_p$ vers une
fonction $g_\infty$, donc $D^{-n} R^{n*}\delta_\nu$ converge vaguement vers
une mesure $\rho_R$. Cette mesure ne d{\'e}pend pas du point $\nu$ choisi.
Notons que pour tout $\nu\in\mathbb{H}_p$ nous pouvons {\'e}crire $\rho_R
= \delta_\nu + \Delta g$ avec $g$ \emph{born{\'e}e}. Ceci implique en particulier
que $\rho_R$ ne charge aucun point de $\mathbb{P}^1(\mathbb{C}_p)$.
Elle ne charge donc pas l'ensemble exceptionnel $\mathcal{E}$ de $R$.
La mesure $\rho_R$ v{\'e}rifie toujours les {\'e}quations d'invariance: $
R^*\rho_R = D\, \rho_R$, et $ R_* \rho_R = \rho_R$.  En
particulier, pour tout $\nu,\nu'\in\mathsf{P}^1(\mathbb{C}_p)$ tels
que $R(\nu') = \nu$ on a $ \rho_R \{\nu'\} =
\frac{\deg(R)}{\deg_R(\nu')}\, \rho_R\{\nu\} $.

On d{\'e}duit directement de ces propri{\'e}t{\'e}s la
Proposition~\ref{prop}.

%On va maintenant d{\'e}montrer la proposition.  Supposons que $R$ ait
%bonne r{\'e}duction.  Alors le point $\nu_0 \in\mathbb{H}_p$ est
%totalement invariant, et donc $R^*\nu_0 = d \nu_0$. Par d{\'e}finition de
%$\rho_R$ on a $\rho_R = \lim d^{-n} R^{n*} \nu = \nu$.  Inversement
%supposons que $\rho_R(\nu) >0$ pour $\nu\in\mathbb{H}_p$.  L'{\'e}quation
%d'invariance implique que $\rho_R(\mu) >0$ pour tout point $\mu$ dans
%$\mathrm{GO}(\nu)$ la grande orbite de $\nu$. On peut donc supposer
%que $\rho_R(\nu) =\max_{\mu \in \mathrm{GO}(\nu)}\rho_R(\mu)$.  Pour
%toute pr{\'e}image de $\nu$, on a alors $\rho_R ( \mu) =
%\frac{\deg(R)}{\deg_R(\mu)} \rho_R(\nu) \le \rho_R(\nu)$ donc
%$\deg_R(\mu) = d$. Le point $\nu$ admet donc une unique pr{\'e}image $\mu$
%par $R$ pour laquelle $\rho_R ( \mu) = \rho_R(\nu)$. Ceci est aussi
%valide pour tout it{\'e}r{\'e} $R^k,k\ge1$. La masse de $\rho_R$ {\'e}tant finie,
%on en d{\'e}duit que $\nu$ est totalement invariant pour un it{\'e}r{\'e} $R^k$ de
%$R$, et donc que $\nu$ lui-m{\^e}me est totalement invariant. La fraction
%rationnelle $R$ admet donc un point de $\mathbb{H}_p$ totalement
%invariant, elle a donc bonne r{\'e}duction par~\cite[]{}.

%

\section{Preuve du th{\'e}or{\`e}me principal}
On a d{\'e}j{\`a} construit la mesure $\rho_R$ et v{\'e}rifi{\'e} qu'elle
{\'e}tait invariante. De l'invariance de $\rho_R$, du fait que son
support est m{\'e}trisable  et qu'elle ne charge pas les points de
$\mathbb{P}^1(\mathbb{C}_p)$, on d{\'e}duit que son support est {\'e}gal
{\`a} l'ensemble de Julia.

Prenons un point $x\in\mathbb{P}^1(\mathbb{C}_p)$ qui n'est pas
totalement invariant par $R$ ou $R^2$, et {\'e}crivons $\delta_x =
\delta_{\nu_0} + \Delta g$ avec $g: \mathbb{H}_p \to \mathbb{R}$. On a
alors, $D^{-n}R^{n*}\delta_x = D^{-n}R^{n*}\delta_{\nu_0} + \Delta (
D^{-n} g \circ R^n)$ pour tout $n\ge0$, donc
$D^{-n}R^{n*}\delta_x\to\rho_R$ d{\`e}s que $D^{-n} g \circ R^n\to 0$
ponctuellement. On va d{\'e}montrer que cette suite de fonctions tends
vers z{\'e}ro. Pour cel{\`a} on va estimer les vitesses de convergence des
points de $\mathbb{H}_p$ vers le bord de l'arbre
$\mathsf{P}^1(\mathbb{C}_p)$ c'est-{\`a}-dire vers
$\mathbb{P}^1(\mathbb{C}_p)$.

La fonction $g$ est support{\'e}e sur le segment reliant $x$ {\`a} $\nu_0$, et
tends vers $+\infty$ en $x$, elle est donc uniform{\'e}ment born{\'e}e
inf{\'e}\-rieurement, et $\liminf D^{-n} g \circ R^n\ge 0$.  Il nous
suffit donc de montrer $\limsup D^{-n} g \circ R^n\le 0$.  Notons
$\mathcal{E}$ l'ensemble des points totalement invariants de $
\mathbb{P}^1(\mathbb{C}_p)$ par $R$ ou $R^2$, et $B(\mathcal{E})$
l'ensemble des points de $ \mathsf{P}^1(\mathbb{C}_p)$ dont les it{\'e}r{\'e}s
convergent vers $\mathcal{E}$.  La fonction $g$ est born{\'e}e dans un
voisinage de $\mathcal{E}$.  Donc $\lim D^{-n} g \circ R^n_*=0$ sur
$B(\mathcal{E})$.

Quitte {\`a} prendre un it{\'e}r{\'e}, on peut supposer que le degr{\'e} local en tout
point hors de $\mathcal{E}$ est strictement inf{\'e}rieur au degr{\'e} de $R$.
Pour chaque point critique $w \in \mathbb{P}^1(\mathbb{C}_p)$ d'ordre
$k(w)\ge 2$, on fixe une petite boule $B(w)\subset
\mathbb{P}^1(\mathbb{C}_p)$, telle que $R (B(w))$ est {\`a} nouveau une
boule, et $R : B(w) \to R(B(w))$ est conjugu{\'e}e {\`a} $z \mapsto z^k$. On
notera $r(w)$ le diam{\`e}tre de $B(w)$. La distance chordale sur
$\mathbb{P}^1(\mathbb{C}_p)$ est ultra-m{\'e}trique, elle se prolonge donc
naturellement en une m{\'e}trique $\mathsf{d}$ sur l'arbre
$\mathsf{P}^1(\mathbb{C}_p)$.  L'enveloppe convexe de $B(w)$ dans
$\mathsf{P}^1(\mathbb{C}_p)$ est encore une boule $\hat{B}(w)$ pour la
m{\'e}trique $\mathsf{d}\,$ de rayon $r(w)$, et pour tout $\mu\in
\hat{B}(w)$, on a $\mathsf{d}\, (R(\mu), \mathbb{P}^1(\mathbb{C}_p))
\ge \mathsf{d}\, (\mu, \mathbb{P}^1(\mathbb{C}_p))^k$.  Hors de la
r{\'e}union des boules $B(w)$ pour $w$ critique, on a par ailleurs
$$\mathsf{d}\, (R(\mu), \mathbb{P}^1(\mathbb{C}_p)) \ge C\,
\mathsf{d}\, (\mu, \mathbb{P}^1(\mathbb{C}_p)) \mbox{ pour une
  constante } C>0~.$$
Hors du bassin d'attraction $B(\mathcal{E})$ des
points exceptionnels, l'orbite d'un point ne peut jamais tomber dans
une boule $\hat{B}(w)$ o{\`u} $k(w)=D$.  On a donc $$\mathsf{d}\,
(R^n_*(\mu), \mathbb{P}^1(\mathbb{C}_p)) \ge \left( C\, \mathsf{d}\,
  (\mu, \mathbb{P}^1(\mathbb{C}_p))\right) ^{(D-1)^n} \mbox{ pour tout
} n~.$$
Par ailleurs, $g(\mu) \le - \log
\mathsf{d}\,(\mu,\mathbb{P}^1(\mathbb{C}_p)) + O(1)$ pour tout
$\mu\in\mathbb{H}_p$, donc $g (R^n_*(\mu)) \lesssim (D-1)^n$ ce qui
prouve $D^{-n} g\circ R^n_*\to0$ hors de $B(\mathcal{E})$.

Finalement, le caract{\`e}re m{\'e}langeant se d{\'e}duit facilement des
propri{\'e}t{\'e}s d'{\'e}qui\-distribution (voir par exemple~\cite[Lemme
4.11]{FG}).

% etc, etc

% Les remerciements sont dans une section, sans num{\`E}rotation

\section*{Remerciements}
% Remerciements - texte ici

Les deux auteurs tiennent {\`a} remercier chaleureusement l'ACI \og dynamique
des applications polynomiales \fg~et tout particuli{\`e}rement Serge Cantat
pour son int{\'e}r{\^e}t et soutien actif \` a ce projet.

\end{document}